\documentclass[12pt,letterpaper]{article}
\usepackage{amsfonts}
\usepackage{sw20lart}
\usepackage{amsthm}
\usepackage{amsmath}

\theoremstyle{definition}
\newtheorem{theorem}{Theorem}
\newtheorem{lemma}[theorem]{Lemma}
\newtheorem{definition}[theorem]{Definition}
\mathcode`\;="203B
\renewenvironment{proof}[1][Proof]{\textbf{#1.} }{\ \rule{0.5em}{0.5em}}
\input tcilatex

\begin{document}

\author{Palle E.T. Jorgensen \\
Mathematics Department, University of Iowa,\\
Iowa City, IA 52242, USA\\
jorgen@math.uiowa.edu\\
http://www.math.uiowa.edu/\symbol{126}jorgen/ \and Anna Paolucci \\
Dipartimento di Scienza dell'Informazione,\\
Universita' di Genova, Italy\\
paolucci@disi.unige.it}
\title{Wavelets in mathematical physics: $q$-oscillators}
\maketitle

\begin{abstract}
We construct representations of a $q$-oscillator algebra by operators on
Fock space on positive matrices. They emerge from a multiresolution scaling
construction used in wavelet analysis. The representations of the Cuntz
Algebra arising from this multiresolution analysis are contained as a
special case in the Fock Space construction.
\end{abstract}

In this paper we establish a connection between multiresolution wavelet
analysis on one hand and representation theory for operator on Hilbert
spaces depending on a real parameter on the other. These operators arise
from a multiresolution wavelet analysis based on Bessel functions. We wish
to develop a framework for the study of creation operators on Hilbert space,
satisfying simple identities, and allowing a Hopf algebra structure.
Examples will include oscillator algebras coming from physical models.

In the first section of the paper, we review the background and the
motivation for the study of the $q$-relations, both as it relates to
problems in mathematics and in physics. On the mathematical side, the
problems concern wavelet analysis and transform theory, especially the
Mellin transform, and on the physics side, they relate to the quon gas of
statistical mechanics. For the construction of the representations, we then
turn to the twisted Fock space and the $q$-oscillator algebra. Our approach
is motivated by wavelet analysis, and it uses a certain loop group. Our main
result is Theorem \ref{ThmOsc.1}.

\section{\label{Int}Introduction}

Some of the results from the papers \cite{jo-bra-ev}, \cite{jo-bra}, \cite%
{jo-kribs} are based on an operator-theoretic approach to wavelet theory
involving representing wavelets in terms of operators in an
infinite-dimensional Hilbert space.

In this paper we introduce an analogous operator approach in a study of a
generalized biorthogonal wavelet, leading to the construction of oscillator
algebras, and more generally Hopf algebras. In \cite{jo-pa} a related class
of representations of the Cuntz algebra $\mathcal{O}_{N}$ have been found
depending on a parameter $q$.

We shall distinguish between two aspects of the study of the $q$-relations:
(i) the $C^{\ast }$-algebra on these relations, and (ii) finding the
representations of them. While the paper \cite{jo-we} by Jorgensen, Schmitt
and Werner covered (i), we shall concentrate here on (ii). This is a
difficult problem: As noted in \cite{jo-we}, the $C^{\ast }$-algebra $%
\mathfrak{A}(q)$ on the $q$-relations is infinite, and its equivalence
classes of irreducible representations do not admit a Borel cross section
for its classification parameters. In \cite{jo-we}, the authors showed that
there is a stability interval $J$ in the $q$-variable such that all the $%
C^{\ast }$-algebra $\mathfrak{A}(q)$ for $q\in J$ are isomorphic to the
Cuntz algebra, see \cite{Bo-Speicher}, and they estimated the size of $J$.

The motivation for the problem (ii) comes from two sources, (a) from
analysis (wavelets, special functions and combinatorics), and (b) from
physics (quantum optics, statistical mechanics, quantum fields and anyons).
We show in this paper how these problems may perhaps be understood better
via the approach of $q$-deformations, and via the study of concrete
mathematical settings where the representations arise naturally.

Other papers which cover representations include \cite{Bo-Speicher}, \cite%
{IS}, and \cite{RS}. The physics of the $q$-relations is outlined in \cite%
{Bo-Speicher}, \cite{pr-vy}, and \cite{RW}. In particular \cite{RW} relates
the $q$-representations to the Gibbs paradox.

The $C^{\ast }$-algebra $\mathcal{O}_{N}$, called the Cuntz $C^{\ast }$%
-algebra on $N$ generators, is universal on the relations 
\begin{equation}
s_{i}^{\ast }s_{j}=\delta _{ij}\mathbf{1},\qquad
\sum_{i=1}^{N}s_{i}s_{i}^{\ast }=\mathbf{1},  \label{eqMR.star}
\end{equation}%
where $\mathbf{1}$ denotes the unit element in $\mathcal{O}_{N}$. If $%
m_{1},\dots ,m_{N}$ are given functions on $\mathbf{T}=\left\{ z\in \mathbf{C%
}:\left| z\right| =1\right\} $, then the operator system 
\begin{equation*}
\left( S_{i}f\right) \left( z\right) =m_{i}\left( z\right) f\left(
z^{N}\right) ,\qquad f\in L^{2}\left( \mathbf{T}\right) ,\;z\in \mathbf{T}%
,\;i=1,\dots ,N
\end{equation*}%
satisfies the Cuntz relations (\ref{eqMR.star}) if and only if the functions
are frequency subband filters for the wavelet multiresolution construction 
\cite{BrJo02b}. Then $m_{1}$ is called the low-pass filter, and the others
filters of the higher frequency bands.

The conditions on the functions may be stated in either one of the following
two equivalent forms (a) or (b): Let $\rho _{N}=e^{i2\pi /N}$.

\renewcommand{\theenumi}{\alph{enumi}} \renewcommand{\labelenumi}{\textup(%
\theenumi)}

(a)\label{MRcondition(1)}The $N\times N$ matrix 
\begin{equation*}
M\left( z\right) =\frac{1}{\sqrt{N}}\left( m_{j}\left( \rho _{N}^{k}z\right)
\right) _{j,k=1}^{N}
\end{equation*}%
is unitary for all $z\in \mathbf{T}$, i.e., 
\begin{equation*}
M\left( z\right) ^{\ast }M\left( z\right) =I_{N},\qquad z\in \mathbf{T}.
\end{equation*}

(b)\label{MRcondition(2)}The $N\times N$ matrix $A\left( z\right) =\left(
A_{j,k}\left( z\right) \right) $ given by 
\begin{equation*}
A_{j,k}\left( z\right) =\frac{1}{N}\sum_{\substack{ w\in \mathbf{T} \\ %
w^{N}=z}}m_{j}\left( w\right) w^{-k}
\end{equation*}%
is unitary for all $z\in \mathbf{T}$, i.e., 
\begin{equation*}
A\left( z\right) ^{\ast }A\left( z\right) =I_{N},\qquad z\in \mathbf{T}.
\end{equation*}

To complete the picture of multiresolution analysis depending on a parameter
we include here the case of the construction of a different multiresolution
via the Mellin transform. We develop a finite scale multiresolution analysis
via Mellin transforms giving rise to wavelets depending on a parameter $q$, $%
0<q<1$.

\section{\label{Mellin}Mellin Transforms}

Let us first recall some facts about Mellin transforms. The Mellin transform
of a function $f\left( x\right) $ is given by 
\begin{equation*}
M\left( f\left( x\right) ;s\right) =\int_{0}^{\infty }x^{s-1}f\left(
x\right) dx.
\end{equation*}%
The inverse transform gives 
\begin{equation*}
f\left( x\right) =M^{-1}\left( F\left( s\right) ;x\right) ,
\end{equation*}%
where we set $M\left( f\left( x\right) ;s\right) =F\left( s\right) $.

The behavior of the Mellin transform under various coordinate transforms in $%
x$ space is given by 
\begin{eqnarray*}
M\left( f\left( ax\right) ;s\right) &=&a^{-s}M\left( f\left( x\right)
;s\right) , \\
M\left( f\left( x^{a}\right) ;s\right) &=&a^{-1}M\left( f\left( x\right) ;%
\frac{s}{a}\right) , \\
M\left( x^{a}f\left( x\right) ;s\right) &=&M\left( f\left( x\right)
;s+a\right) ,
\end{eqnarray*}%
as can be easily checked.

Let us give some preliminaries on standard multiresolution wavelet analysis
of scale $N$. Define scaling by $N$ on $L^2\left( \mathbf{R}\right) $ by 
\begin{equation*}
U\left( \xi \right) \left( x\right) =\xi \left( N^{-1}x\right)
\end{equation*}
and translation by $1$ by 
\begin{equation*}
T\left( \xi \right) \left( x\right) =\xi \left( x-1\right) .
\end{equation*}
A scaling function is a function $\phi \in L^2\left( \mathbf{R}\right) $
such that if $V_0$ is the closed linear span of all translated $T^k\phi $, $%
k\in \mathbf{Z}$, then $\phi $ has the following properties:

\begin{enumerate}
\item \label{scalfunc(1)} $\left\{ T^k\phi :k\in \mathbf{Z}\right\} $ is an
orthonormal set in $V_0$,

\item \label{scalfunc(2)} $U\phi \in V_0$,

\item \label{scalfunc(3)} $\bigwedge_{n\in \mathbf{Z}}U^nV_0=\left\{
0\right\} $,

\item \label{scalfunc(4)} $\bigvee_{n\in \mathbf{Z}}U^nV_0=L^2\left( \mathbf{%
R}\right) $,
\end{enumerate}

\noindent where the symbol $\bigvee $ means ``closed linear span'', and $%
\bigwedge $ ``intersection of subspaces''. We shall further use the notation 
$U_{j}:=U^{j}V_{0}$, $j\in \mathbf{Z}$. The system \ref{scalfunc(1)}.--\ref%
{scalfunc(4)}.\ is called a multiresolution analysis (MRA).

We then have the following result

\begin{theorem}
\label{ThmInt.1} Let $\Phi \in L^{2}\left( \mathbf{R}\right) $ be the Haar
function, $0<q<1$, $q\in \mathbf{R}$. There exists a sequence of subspaces $%
U_{j}$, $j\in \mathbf{Z}$, such that 
\begin{equation}
\cdots U_{j}\subset U_{j-1}\subset \dots \subset U_{0}\subset U_{-1}\subset
\cdots  \label{EqSub}
\end{equation}%
giving a multiresolution satisfying the above properties.
\end{theorem}

\begin{proof}
To show that the system in (\ref{EqSub}) is a multi-resolution, we need to
establish a scaling operator which moves in steps of one through the ladder
of resolution spaces $U_{j}$; we must identify the Hilbert space $\mathcal{H}
$; and finally we must show that $\bigvee_{j}U_{j}=\mathcal{H}$, and $%
\bigwedge_{j}U_{j}=\{0\}$ where $\bigvee $ and $\bigwedge $ are the usual
lattice operations in Hilbert space.

Let $\Phi \in L^{2}\left( \mathbf{R}\right) $ be the Haar function, $0<q<1$, 
$x\in \mathbf{R}$. We want to prove that there exists a sequence of
subspaces $U_{j}$ giving rise to a multiresolution satisfying the above
properties. Define 
\begin{equation*}
V_{0}=\bigvee \left\{ \Phi \left( q^{k}x\right) :k\in \mathbf{Z}%
,\;q<x<1\right\} .
\end{equation*}%
Let $U\xi \left( x\right) =\frac{1}{\sqrt{N}_{\mathstrut }}\xi \left(
x^{N}\right) $ be the scaling operator and $T\xi \left( x\right) =\xi \left(
qx\right) $ playing the role of the ``translation operator''. In this case $%
T $ is a dilation. Then the set $\left\{ \Phi \left( q^{k}x\right) :k\in 
\mathbf{Z}\right\} $ is orthonormal in $L^{2}\left( \mathbf{R}_{+}\right) $, 
$q\in \mathbf{R}$, $0<q<1$, since the intervals are disjoint for $h\neq k$, $%
h,k\in \mathbf{Z}$. It is easy to verify that $U\Phi \in V_{0}$, 
\begin{equation}
U\Phi \left( x\right) =\sum_{k}a_{k}\Phi \left( q^{k}x\right) .
\label{uno(1)}
\end{equation}%
Then $\bigwedge_{n\in \mathbf{Z}}U^{n}V_{0}=\left\{ 1\right\} $. Suppose $%
\xi \in U^{n}V_{0}$, with $\xi \left( x\right) =U^{n}\Phi \left( x\right)
=\sum_{k}a_{k}\Phi \left( \left( q^{k}x\right) ^{Nn}\right) \frac{1}{Nn}$.
Since any $\xi \in L^{2}\left( \mathbf{R}_{+}\right) $ can be approximated
by step functions and $\mathbf{R}_{+}=\bigcup_{k,n}\left[ q^{\left(
k+1\right) Nn},q^{kNn}\right] $, it follows that 
\begin{equation*}
\bigvee_{n\in \mathbf{Z}}U^{n}V_{0}=L^{2}\left( \mathbf{R}_{+}\right) .
\end{equation*}%
Hence the above properties are satisfied and the sequence $V_{j}=U^{j}V_{0}$
of subspaces $V_{j}$ associated to $\Phi $ defines a multiresolution. Let $%
\Phi \in L^{2}\left( \mathbf{R}\right) $ satisfy (\ref{uno(1)}). Define the
operator $W_{\Phi }\colon L^{2}\left( \mathbf{R}\right) \longrightarrow
L^{2}\left( \mathbf{R}\right) $, called the wave operator, by 
\begin{equation*}
\left( W_{\Phi }f\right) \left( x\right) =\sum_{k\in \mathbf{Z}}c_{k}\Phi
\left( q^{k}x\right) =\sum_{k\in \mathbf{Z}}M\left( f\left( q^{k}x\right)
;s\right) \Phi \left( q^{k}x\right) .
\end{equation*}

\begin{lemma}
\label{LemInt.2}The operator $W_\Phi $ satisfies 
\begin{equation}
UW_\Phi =W_\Phi S ,  \label{tre}
\end{equation}
where $S$ is given by 
\begin{equation*}
M\left( Sf\right) \left( s\right) = \sum_{j}a_{j}q^{-sj}M\left( f;\frac
sN\right) .
\end{equation*}
\end{lemma}

\begin{proof}
First, we have 
\begin{eqnarray*}
\left( UW_\Phi f\right) \left( x\right) &=&\frac 1{\sqrt{N}}\left( W_\Phi
f\right) \left( x^N\right) =\frac 1{\sqrt{N}}\sum_{k\in \mathbf{Z}}M\left(
f\left( q^{kN}x^N\right) ;s\right) \Phi \left( q^{kN}x^N\right) \\
\ &=&\frac 1{\sqrt{N}}\sum_{k\in \mathbf{Z}}q^{-ksN}M\left( f\left( x\right)
;\frac sN\right) \Phi \left( q^{kN}x^N\right) \\
\ &=&\frac 1{\sqrt{N}}\sum_{k\in \mathbf{Z}}\sum_{j\in \mathbf{Z}%
}q^{-ksN-sj}a_jM\left( f\left( x\right) ;\frac sN\right) \Phi \left(
q^{kN+j}x^N\right) \\
\ &=&\frac 1{\sqrt{N}}\sum_{k\in \mathbf{Z}}\left[ \sum_{l\in \mathbf{Z}%
}a_{l-Nk}q^{-ls}M\left( f\left( x\right) ;s\right) \right] \Phi \left(
q^lx^N\right) .
\end{eqnarray*}
On the other hand 
\begin{eqnarray*}
W_\Phi \left( Sf\right) \left( z\right) &=&\sum_{k\in \mathbf{Z}}M\left(
\left( Sf\right) \left( xq^k\right) ;s\right) \Phi \left( q^kx\right) \\
&=&\sum_{k\in \mathbf{Z}}M\left( \left( Sf\right) \left( xq^k\right)
;s\right) \Phi \left( q^{kN+j}x^N\right) \\
\ &=&\sum_{k\in \mathbf{Z}}\sum_{j\in \mathbf{Z}}q^{-kNs-sj}a_jM\left(
f\left( x\right) ;\frac sN\right) \Phi \left( q^{kN+j}x^N\right) \\
\ &=&\sum_{k\in \mathbf{Z}}\sum_{l\in \mathbf{Z}}a_{l-kN}q^{-ls}M\left(
f\left( x\right) ;\frac sN\right) \Phi \left( q^lx^N\right) .
\end{eqnarray*}
Thus 
\begin{equation*}
UW_\Phi =W_\Phi S
\end{equation*}
which is the identity (\ref{tre}), i.e., $W_\Phi $ intertwines the operators 
$U$ and $S$.
\end{proof}

Using the fact that $\left\{ T^k\Phi :k\in \mathbf{Z}\right\} $ is an
orthonormal set in $L^2\left( \mathbf{R}\right) $ we may define 
\begin{equation*}
F_\Phi\colon V_0\longrightarrow L^2\left( \mathbf{R}\right) ,
\end{equation*}
an isometry which sends 
\begin{equation*}
\xi \longmapsto m\left( s\right) =\sum_ka_kq^{-ks}, \qquad \xi \left(
s\right) =\sum_ka_k\Phi \left( q^{-k}x\right) .
\end{equation*}
Then 
\begin{equation*}
M\left( \xi \left( x\right) ;s\right) =m\left( s\right) M\left( \Phi \left(
x\right) ;s\right)
\end{equation*}
where $M\left( \xi \left( x\right) ;s\right) $ is the Mellin transform.

If $\xi \in V_{-1}=U^{-1}V_0$ then $U\xi \in V_0$ so we define 
\begin{equation*}
m_\xi =F_\Phi \left( U\xi \right) \in L^2\left( \mathbf{R}\right)
\end{equation*}
and then 
\begin{eqnarray*}
M\left( U\xi \left( x\right) ;s\right) &=&M\left( \sum_ka_k\Phi \left(
\left( q^kx\right) ^N\right) ;s\right) =\sum_ka_kM\left( \Phi \left( \left(
q^kx\right) ^N\right) ;s\right) \\
&=&\sum_ka_kM\left( \Phi \left( q^{kN}x^N\right) ;s\right) =\sum_k\frac
1Na_kq^{-ksN}M\left( \Phi \left( x\right) ;\frac sN\right) .
\end{eqnarray*}
Thus 
\begin{equation}
NM\left( U\xi \left( x^N\right) ;s\right) =m_\xi \left( s\right) M\left(
\Phi \left( x\right) ;s\right) .  \label{due(1)}
\end{equation}

Observe that in this case we don't have unitarity of the matrix $A$
generated by the filters $m_\xi $. Thus we don't have representations of the
Cuntz algebra. This brings out the question of what are the natural operator
relations associated to this MRA. One particular choice is given by the
oscillator algebras.

For any MRA wavelet construction as found above we construct an appropriate
oscillator algebra such that the filters make up the eigenfunctions system
of the oscillator hamiltonian, i.e., the energy operator.

In fact the non-unitarity of the matrix given by the MRA filters provides
the eigenfunctions for the energy operator. The $q$-oscillator algebra is
given by the operators $a^{-}$, $a^{+}$ and the number operator. By a
standard argument it is possible to build up the following selfadjoint
operators: 
\begin{equation*}
X =\frac{\sqrt{2}}2\left( a^{+}+a^{-}\right) \text{\quad and \quad }P =\frac{%
\sqrt{2}}2\left( a^{-}-a^{+}\right)
\end{equation*}
which are the momentum and the position operators. Given the above MRA there
exists a family of $q$-oscillators that can be represented via the filters.

Let $H=m_0^2(s) + m_0^2(\sigma(s))$ be the energy operator. Let 
\begin{equation*}
a^{-}(s)a^{+}(s)=m_0(\sigma(s))^2
\end{equation*}
and 
\begin{equation*}
a^{+}(s)a^{-}(s)=m_0(s)^2.
\end{equation*}

Then it follows that the algebra is generated by $a^{+}$, $a^{-}$, $N$ and
the generators satisfy the following relations 
\begin{equation*}
\left[ a^{-}(s),a^{+}(s)\right] =f(N+1) - f(N)= f(N) [s]_{q^{-2}}
\end{equation*}

\begin{equation*}
f(N) a^{-}(s) = a^{-}(s) f(N-1),
\end{equation*}

\begin{equation*}
f(N) a^{+}(s) = a^{+}(s) f(N+1),
\end{equation*}
where $f(N+1)-f(N) = \sum_k b_k q^{-2ks}[s]_{q^{-2}}$.

Thus we have 
\begin{equation*}
\left[ a^{-}(s),a^{+}(s)\right] =m_0(s)^2\left[ s\right] _{q^{-2}}
\end{equation*}
where $\left[ s\right] _q=\frac{1-q^s}{1-q}$. Then using MRA we can
construct a family of operators $\left\{ a_i^{-}(s),a_i^{+}(s)\right\} $ for 
$i=1,\dots ,N$ satisfying: 
\begin{equation*}
\left[ a_i^{-}(s),a_i^{+}(s)\right] =m_i(s)^2\left[ s\right] _{q^{-2}}
\end{equation*}

A representation of the operators $a^{-}(s)$ and $a^{+}(s)$ can be realized
in a Fock space as follows: 
\begin{equation*}
a^{-}(s)\left| e_k\right\rangle =\frac q{1-q}\sum_kb_kq^{-ks}\left|
e_{k-1}\right\rangle
\end{equation*}
and: 
\begin{equation*}
a^{+}(s)\left| e_k\right\rangle =\frac q{1+q}\sum b_kq^{-\left( k+1\right)
s}\left| e_{k+1}\right\rangle.
\end{equation*}
\end{proof}

\section{\label{Twi}Twisted Fock Space}

In this Section we introduce a new Fock space construction, see \cite%
{jo-kribs}, which may provide the appropriate framework for studying wavelet
representations of certain $q$-oscillator algebras.

\begin{definition}
\label{DefTwi.1} The full Fock space over $\mathbf{C}^{N}$ where $N$ is a
fixed positive integer with $N\geq 2$, is the orthogonal direct sum of
Hilbert spaces given by 
\begin{equation*}
K=\left( \;\sideset{}{^{\smash{\oplus}}}{\sum}\limits_{k=-\infty
}^{-1}\left( \mathbf{C}^{N}\right) ^{\otimes -k}\right) \oplus \mathbf{C}%
=\cdots \oplus \left( \mathbf{C}^{N}\otimes \mathbf{C}^{N}\right) \oplus
\left( \mathbf{C}^{N}\right) \oplus \mathbf{C}.
\end{equation*}%
The term $\mathbf{C}$ in the summand on the right designates the vacuum
vector (in the formula we omit a special symbol $\Omega $ for the vacuum
vector). Let $\left\{ \xi _{1},\dots ,\xi _{N}\right\} $ be a fixed
orthonormal basis for $\mathbf{C}^{N}$. Then $K$ is an infinite-dimensional
Hilbert space with orthonormal basis given by $\left\{ \xi _{i_{1}},\dots
,\xi _{i_{k}}:1\leq i_{1},\dots ,i_{k}\leq N,\;k\geq 1\right\} \cup \left\{
\Omega \right\} $.
\end{definition}

\subsection{Construction of Twisted Fock space}

Let $K$ be a Hilbert space. We take the tensor product of the full Fock
space with $H$, then we define a ``new'' inner product $\langle \, \cdot \,
,\, \cdot \, \rangle _\Phi $ by using a completely positive map from the
complex matrices into $B\left( H\right) $ (i.e., a positive matrix with
entries in $B\left( H\right) $). By following \cite{jo-kribs} we construct
the twisted Fock space as follows. Let $\Phi\colon M_N\longrightarrow
B\left( H\right) $ be the completely positive map which we will define
later. We define the $N$-variable pre-Fock space over $K$ to be the vector
space of finite sums 
\begin{equation*}
T_N\left( H\right) =\left\{ \sum_{\left| w\right| \le k}w\otimes h_w:w\in 
\mathbf{F}_N^{+},\;k\ge 1,\;h_w\in H\right\}
\end{equation*}
where $\mathbf{F}_N^{+}$ is the unital free-semigroup on $N$ non-commuting
letters $\left\{ 1,2,\dots ,N\right\} $ with unit $e$. We can think of the
full Fock space as $l^2\left( \mathbf{F}_N^{+}\right) $ where an orthonormal
basis is given by the vectors $\left\{ \xi _w:w\in \mathbf{F}_N^{+}\right\} $
corresponding to the words $w$, $\left| w\right| $ is the word of length
zero or empty word, $\xi _w=\xi _{i_1}\otimes \dots \otimes \xi _{i_n}$, $%
w=i_1\cdots i_k\in \mathbf{F}_N^{+}$. Then a vector $\left( i_1,\dots
,i_k\right) \otimes h$ with $w=i_1\cdots i_k\in \mathbf{F}_N^{+}$
corresponds to the vector $\xi _{i_1}\otimes \dots \otimes \xi _{i_n}\otimes
h$ in $\left( \mathbf{C}^N\right) ^{\otimes k}\otimes H$.

Let $\Phi $ be the completely positive map $\Phi\colon M_N\longrightarrow
B\left( H\right) $. Define a form 
\begin{equation*}
\langle \, \cdot \, ,\, \cdot \, \rangle _\Phi\colon T_N\left( H\right)
\times T_N\left( H\right) \longrightarrow \mathbf{C}
\end{equation*}
as follows. For $w,w^{\prime }\in \mathbf{F}_N^{+}$, $h,h^{\prime }\in H$

\begin{itemize}
\item[i)] $\langle e\otimes h,e^{\prime }\otimes h^{\prime }\rangle _\Phi
=\langle h\mid h^{\prime }\rangle $;

\item[ii)] if $\left| w\right| \ne \left| w^{\prime }\right| $ then $\langle
w\otimes h,w^{\prime }\otimes h^{\prime }\rangle _\Phi =0$;

\item[iii)] if $w=i_1\cdots i_k$, $w^{\prime }=i_1^{\prime }\cdots
i_k^{\prime }$, then 
\begin{equation*}
\langle w\otimes h,w^{\prime }\otimes h^{\prime }\rangle _\Phi =\left\langle
h\mid \Phi \left( e_{i_1i_1^{\prime }}\otimes \dots \otimes
e_{i_ki_k^{\prime }}\right) h^{\prime }\right\rangle .
\end{equation*}
\end{itemize}

Extend $\langle \, \cdot \, ,\, \cdot \, \rangle _\Phi $ to $T_N\left(
H\right) \times T_N\left( H\right) $ as linear in the first variable and
conjugate linear in the second one. From Theorem 4.5 of \cite{jo-kribs} the
form $\langle \, \cdot \, ,\, \cdot \, \rangle _\Phi $ is positive
semi-definite on $T_N\left( H\right) $.

\begin{definition}
Let $N_\Phi =\left\{ x\in T_N\left( H\right) :\langle x\mid x\rangle _\Phi
=0\right\} $ be the kernel of the form $\langle \,\cdot \,,\,\cdot \,\rangle
_\Phi $. Define the Fock space of $\Phi $ over $H$ to be the Hilbert space
completion 
\begin{equation*}
F_N\left( K,\Phi \right) =\overline{T_N\left( H\right) /N_\Phi }^{\langle
\,\cdot \,,\,\cdot \,\rangle _\Phi }.
\end{equation*}
The left creation operators $T=\left( T_1,\dots ,T_N\right) $ on $\mathcal{F}%
_N\left( H,\Phi \right) $ are linear transformations defined by 
\begin{equation*}
T_i\left( w\otimes h+N_\Phi \right) =\left( iw\right) \otimes h+N_\Phi .
\end{equation*}
These operators are well-defined and $T_i\left( N_\Phi \right) \subset
N_\Phi $, $1\le i\le N$.

Let $S_i$ be the $i$-th creation operator on the first space defined as
follows 
\begin{equation*}
S_ix=S_i\left( \eta _{i_1}\otimes \dots \otimes \eta _{i_k}\otimes h\right)
=\eta _i\otimes \eta _{i_1}\otimes \dots \otimes \eta _{i_k}\otimes h
\end{equation*}
for vectors $x=\eta _{i_1}\otimes \dots \otimes \eta _{i_k}\otimes h$.
Define a twisted new Fock space as 
\begin{equation*}
T_N\left( H\right) =\left\{ \sum_{\left| w\right| \le k}w\otimes h_w:w\in 
\mathbf{F}_N^{+},\;k\ge 1,\;h_w\in H\right\}
\end{equation*}
as before.

We define a form 
\begin{equation*}
\langle \,\cdot \,,\,\cdot \,\rangle _\Phi \colon T_N\left( H\right) \times
T_N\left( H\right) \longrightarrow \mathbf{C}
\end{equation*}
as follows. For $w,w^{\prime }\in \mathbf{F}_N^{+}$, $h,h^{\prime }\in H$

\begin{itemize}
\item[i)] $\langle e\otimes h,e^{\prime }\otimes h^{\prime }\rangle _\Phi
=\langle h\mid h^{\prime }\rangle $;

\item[ii)] if $\left| w\right| \ne \left| w^{\prime }\right| $ then $\langle
w\otimes h,w^{\prime }\otimes h^{\prime }\rangle _\Phi =0$;

\item[iii)] if $w=i_1\cdots i_k$, $w^{\prime }=i_1^{\prime }\cdots
i_k^{\prime }$, then 
\begin{eqnarray*}
& & \langle w\otimes h,w^{\prime }\otimes h^{\prime }\rangle _\Phi \\
& & \qquad =\left\langle h\bigm| \Phi \left( e_{i_1\sigma \left( i_1^{\prime
}\right) }\otimes \dots \otimes e_{i_k\sigma \left( i_k^{\prime }\right)
}\right) h^{\prime }\right\rangle ,
\end{eqnarray*}
with $i\left( \sigma \right) =\#\left\{ \left( i,j\right) \in \left\{
1,\dots ,N\right\} ^2:i<j,\;\sigma \left( i\right) >\sigma \left( j\right)
\right\} $.
\end{itemize}

We need the map $\phi $ being completely positive and then $\langle \,\cdot
\,,\,\cdot \,\rangle _\Phi $ is positive semi-definite: 
\begin{equation*}
\tilde \Phi \left( e_i\otimes e_j\right) =\Phi \left( e_i\otimes P_q^{\left(
N\right) }e_j\right)
\end{equation*}
\end{definition}

Let us now turn to the construction of the multiresolution wavelet analysis.

Suppose we have two pairs of filters and then we have two pairs of scaling
functions plus wavelet $\Phi $, $\Psi $ and $\tilde \Phi $, $\tilde \Psi $.
They are defined by 
\begin{equation*}
\begin{array}{rclrcl}
\hat \Phi \left( \xi \right) & = & m_0\left( \xi /N\right) \hat \Phi \left(
\xi /N\right) , & \qquad \hat \Psi \left( \xi \right) & = & m_1\left( \xi
/N\right) \hat \Phi \left( \xi /N\right) _{\mathstrut}, \\ 
\Hat{\Tilde{\Phi}}\left( \xi \right) & = & \tilde m_0\left( \xi /N\right) 
\Hat{\Tilde{\Phi}}\left( \xi /N\right) , & \qquad \Hat{\Tilde{\Psi}} \left(
\xi \right) & = & \tilde m_1\left( \xi /N\right) \Hat{\Tilde{\Phi}} \left(
\xi /N\right) ^{\mathstrut}.%
\end{array}%
\end{equation*}
We want to take the direct sum of two MRA. We have the following two
sequences of successive approximation spaces $U_j$ and $\tilde U_j$. The
closed subspaces satisfy 
\begin{equation}
\cdots V_2\subset V_1\subset V_0\subset V_{-1}\subset V_{-2}\subset \cdots
\label{uno(2)}
\end{equation}
with 
\begin{equation*}
\bigvee _{j\in \mathbf{Z}}V_j=L^2\left( \mathbf{C}\right) ,\qquad\bigwedge
_{j\in \mathbf{Z}}V_j=\left\{ 0\right\} ,
\end{equation*}
\begin{equation}
\cdots \tilde V_2\subset \tilde V_1\subset \tilde V_0\subset \tilde
V_{-1}\subset \tilde V_{-2}\subset \cdots  \label{due(2)}
\end{equation}
with 
\begin{equation*}
\bigvee _{j\in \mathbf{Z}}\tilde V_j=L^2\left( \mathbf{C}\right)
,\qquad\bigwedge _{j\in \mathbf{Z}}\tilde V_j=\left\{ 0\right\} .
\end{equation*}
Formulae (\ref{uno(2)}) and (\ref{due(2)}) have the additional requirements 
\begin{equation*}
f\in V_j\;\Leftrightarrow \;f\left( N^j \, \cdot \, \right) \in
V_0,\;\;\;\tilde f\in \tilde V_j\;\Leftrightarrow \;\tilde f\left( N^j \,
\cdot \, \right) \in \tilde V_0,
\end{equation*}
i.e., all spaces are scaled versions of the central space $V_0$ and $\tilde
V_0$ respectively. For every $j\in \mathbf{Z}$, define $W_j$ to be the
orthogonal complement of $V_j$ in $V_{j-1}$ and $\tilde{W}_j$ the orthogonal
complement of $\tilde V_j$ in $\tilde V_{j-1}$. We have 
\begin{equation*}
V_{j-1}=V_j\oplus W_j,\;\;\tilde V_{j-1}=\tilde V_j\oplus \tilde W_j.
\end{equation*}
Also $W_j\perp W_{j^{\prime }}$ if $j\ne j^{\prime }$ and $\tilde W_j\perp
\tilde W_{j^{\prime }}$ if $j\ne j^{\prime }$. Define a sequence of
successive approximation spaces 
\begin{equation*}
Y_1=V_1\oplus \tilde V_1,\;\;Y_0=V_0\oplus \tilde
V_0,\;\;\;Y_{-1}=V_{-1}\oplus \tilde V_{-1},\;\;\;Y_{-2}=V_{-2}\oplus \tilde
V_{-2},\;\dots
\end{equation*}
such that if $f\left( t,s\right) $ is in $Y_j$ then $f\left( Nt,Ns\right) $
and all $f\left( t-k,s-k\right) $ are in $Y_{j+1}$. Define $\hat W_j$ to be
the orthogonal complement of $V_j\cap \tilde V_j$ in $Y_{j-1}$. Then 
\begin{equation*}
\hat W_j=\left( V_j\cap \tilde V_j\right) ^{\perp }\cap Y_{j-1}=V_j^{\perp
}\oplus \tilde V_j^{\perp }\cap Y_{j-1}.
\end{equation*}
Thus 
\begin{eqnarray*}
Y_j &=&\left( V_j\oplus \tilde V_j\right) =\left( V_{j-1}\oplus W_j\right)
\oplus \left( \tilde V_{j-1}\oplus \tilde W_j\right) \\
\ &=&\left( V_{j-1}\oplus \tilde V_{j-1}\right) \oplus \left( W_j\oplus
\tilde W_j\right) =Y_{j-1}\oplus \hat W_j
\end{eqnarray*}
which implies $Y_j=Y_{j-1}\oplus \hat W_j$, where we set $\hat W_j=W_j\oplus
\tilde W_j$. Then 
\begin{eqnarray*}
\hat W_j &=&W_j\oplus \hat W_j=\left( V_j^{\perp }\cap V_{j-1}\right) \oplus
\left( \tilde V_j^{\perp }\cap \tilde V_{j-1}\right) \\
\ &=&\left( V_j^{\perp }\oplus \tilde V_j^{\perp }\right) \cap \left(
V_{j-1}\oplus \tilde V_{j-1}\right) =\left( V_j\cap \tilde V_j\right)
^{\perp }\cap Y_{j-1},
\end{eqnarray*}
so that $\hat W_j$ is the orthogonal complement of $V_j\cap \hat V_j$ in $%
Y_{j-1}$.

Thus $L^2\left( \mathbf{C}^2\right) =\bigoplus_{j\in \mathbf{Z}}\hat W_j$.
The basic point of MRA is that whenever a collection of closed subspaces
satisfies 
\begin{equation*}
\cdots \subset Y_2\subset Y_1\subset Y_0\subset Y_{-1}\subset Y_{-2}\subset
\cdots
\end{equation*}
\begin{equation*}
\bigvee _{j\in \mathbf{Z}}Y_j=L^2\left( \mathbf{C}^2\right) ,\qquad\bigwedge
_{j\in \mathbf{Z}}Y_j=\left\{ 0\right\} ,
\end{equation*}
\begin{equation*}
f\in Y_j\text{ if and only if }f(N^j \, \cdot \, )\in Y_0
\end{equation*}
\begin{equation*}
f\in Y_0\text{ then }( \, \cdot \, -k)\in Y_0\text{ for every }k\in \mathbf{Z%
}
\end{equation*}
\begin{eqnarray*}
& &\text{there exists }\phi \in Y_{0}\text{ such that }\{\phi _{0,k}:k\in 
\mathbf{Z}\}\text{ is an orthonormal basis in }Y_0 \\
& &\text{where }\phi _{j,k}(z) =N^{-j}\phi (N^{-j}z-k)
\end{eqnarray*}

Then there exists an orthonormal wavelet basis $\{\psi _{j,k}:j,k\in \mathbf{%
Z}\}$ of $L^{2}(\mathbf{C}^{2})$, and we may therefore define the following
isomorphism $\eta \colon \mathbf{C}^{2}\rightarrow \mathcal{H}$, where $%
\mathcal{H}$ denotes the quaternions, by $\eta (z_{1},z_{2})=z_{1}+z_{2}e_{2}
$ and $e_{2}=(0,0,1,0)\in \mathbf{R}^{4}$. Consider now the space $L^{2}(%
\mathcal{H})$ equipped with the usual quaternion inner product: Given the
filters $m_{i}$ and $\tilde{m}_{j}$ we define the following two matrix
functions $A$ and $\tilde{A}$ by 
\begin{equation*}
A_{k,l}(z)=\frac{1}{N}\sum_{w^{N}=z}w^{-l}m_{k}(w)
\end{equation*}%
and 
\begin{equation*}
\tilde{A}_{k,l}(z)=\frac{1}{N}\sum_{w^{N}=z}w^{-l}\tilde{m}_{k}(w)
\end{equation*}%
respectively.

They satisfy the following biorthogonality conditions 
\begin{equation*}
\sum_{k=0}^{N-1}\overline{A_{k,i}(z)}\,\tilde{A}_{k,j}(z)=\delta _{i,j}
\end{equation*}
and 
\begin{eqnarray*}
\frac 1N\sum_{w^N=z}\overline{m_i(w)}\,m_j(w) &=&\delta _{i,j} , \\
\frac 1N\sum_{w^N=z}\overline{\tilde{m}_i(w)}\,\tilde{m}_j(w) . &=&\delta
_{i,j}
\end{eqnarray*}

Take $B=A\oplus A^{*}$ and $\tilde B=\tilde A\oplus \tilde A^{*}$ the matrix
functions associated to the MRA's with filters $m_i$ and $\tilde m_i$
respectively.

It can be easily checked for $N=2$ that the following equation is satisfied: 
\begin{equation*}
B\tilde B^{*}+\tilde B^{*}B=\mathbf{1}
\end{equation*}

\begin{theorem}
\label{ThmCon.2} Let $S=(S_0,S_1,\dots ,S_N)$ and $\tilde S=(\tilde
S_0,\tilde S_1,\dots ,\tilde S_N)$ be a pair of wavelet representations on $%
H=L^2(\mathbf{C)}$ with invertible loop matrices $A$ and $\tilde A$
respectively. Let $\mathbf{S}=\left( 
\begin{array}{cc}
S & \tilde S%
\end{array}
\right) $ be the matrix associated to $S$ and $\tilde S$ and let $P=\mathbf{S%
}^{*}\mathbf{S}+\mathbf{S}\mathbf{S}^{*}$.

Let $T=(T_1,T_2,\dots ,T_{N-1},\tilde T_1,\tilde T_2,\dots ,\tilde T_{N-1})$
be the creation operator on $\mathcal{F}_{2N}\left( H,P\right) $. Then:

$T_i^{*}T_j|_H=\left( S_{i-1}^{*}S_{j-1}\right) =\left( AA^{*}\right) _{i,j}$

$\tilde T_i^{*}\tilde T_j|_H=\left( \tilde S_{i-1}^{*}\tilde S_{j-1}\right)
=\left( \tilde A\tilde A^{*}\right) _{i,j}$

$T_i\tilde T_j^{*}|_H+\tilde T_j^{*}T_i|_H=\delta _{i,j}\mathbf{1}$ where $H$
denotes the space $L^2(\mathcal{H})$.

Hence the $*$-algebra generated by $\tilde T_j$, $j=1,\dots ,N$ is an
oscillator algebra. It has a representation containing the Cuntz-Toeplitz
isometries.
\end{theorem}

\begin{proof}
It is easy to check that the fermion algebra relations hold for the direct
sum of wavelet representations.
\end{proof}

\subsection{Construction of $q$-oscillator algebras}

Let us consider the system of wavelet representations $S= (S_0,S_1,\dots
,S_N)$ and $\tilde{S}=(\tilde{S}_0,\tilde{S}_1,\dots ,\tilde{S}_N)$
depending on a real parameter $q$ on $L^2(\mathbf{C)}$ with invertible loop
matrices $A$ and $\tilde{A}$ respectively.

Let $S=(S_0,S_1,\dots ,S_N)$ and $\tilde S=(\tilde S_0,\tilde S_1,\dots
,\tilde S_N)$ be a pair of wavelet representations on $L^2(\mathbf{C)}$ with
invertible loop matrices $A$ and $\tilde A$ respectively. Let us assume that
one of them, say $\tilde S=(\tilde S_0,\tilde S_1,\dots ,\tilde S_N)$,
depends on a real parameter $q$ and the matrix $A$ is given by: 
\begin{equation*}
\tilde A_{k,l}(z)=\sum_{w^N=z}\left( qw\right) ^{-l}\tilde m_k(w)
\end{equation*}
where the $\tilde m_j$ are the filters of the MRA depending on a parameter $%
q $ as constructed in \cite{jo-pa}. Then $\tilde A_{k,l}$ is a unitary
matrix function. The unitarity of the $\left\{ \left( 1-q^{2N}\right) \tilde
m_j\left( tq^j\right) \right\} _{i,j=0,\dots ,N-1}$ is equivalent to 
\begin{equation*}
\sum_{k}\tilde A_{i,k}(z)\,\overline{\tilde A_{j,k}(z)}=\left(
1-q^{2N}\right) \sum_{w^N=z}\tilde m_i(w)\,\overline{\tilde m_j(w)}.
\end{equation*}

We assume then a generalized biorthogonality holds 
\begin{equation*}
\sum_{k}A_{i,k}(z)\,\overline{\tilde A_{j,k}(z)}=\left( 1-q^N\right)
\sum_{w^N=z}m_i(w)\,\overline{\tilde m_j(w)}=\delta _{i,j}1.
\end{equation*}
Let $\mathbf{S}$ be the matrix defined as above. Let $P=\mathbf{S}^{*}%
\mathbf{S}$ be the matrix $2N\times 2N$ in $B(\mathcal{H})$ determined by
the wavelet representations such that 
\begin{equation*}
\mathbf{S}^{*}\mathbf{S}=\left( 
\begin{array}{cc}
S^{*}S & 0 \\ 
0 & \tilde S^{*}\tilde S%
\end{array}
\right)
\end{equation*}

Then we have the following theorem.

\begin{theorem}
\label{ThmOsc.1} Let $S=(S_0,S_1,\dots ,S_N)$ and $\tilde S=(\tilde
S_0,\tilde S_1,\dots ,\tilde S_N)$ be a pair of wavelet representations on $%
L^2(\mathbf{C)}$ with invertible loop matrices $A$ and $\tilde A$
respectively. Let us assume that one of them, say $\tilde S=(\tilde
S_0,\tilde S_1,\dots ,\tilde S_N)$, depends on a real parameter $q$ and the
matrix $A$ is given by: 
\begin{equation*}
\tilde A_{k,l}(z)=\sum_{w^N=z}\left( qw\right) ^{-l}\tilde m_k(w)
\end{equation*}

Let $\mathbf{S}$ be as above and let $P=\mathbf{S}^{*}\mathbf{S}$ be the
matrix $2N\times 2N$ in $B(\mathcal{H})$.

Let $T=(T_1,T_2,\dots ,T_{N-1},\tilde T_1,\tilde T_2,\dots ,\tilde T_{N-1})$
be the creation operator in $\mathcal{F}_{2N}\left( H,P\right) $. Then we
have:

$T_i^{*}T_j|_H=\left( S_{i-1}^{*}S_{j-1}\right) =\left( AA^{*}\right) _{i,j}$

$\tilde T_i^{*}\tilde T_j|_H=\left( \tilde S_{i-1}^{*}\tilde S_{j-1}\right)
=\left( \tilde A\tilde A^{*}\right) _{i,j}$

$\tilde T_iT_j^{*}|_H-T_j^{*}\tilde T_i|_H=\delta _{i,j}\left[ N\right] _q%
\mathbf{1}$
\end{theorem}

\begin{proof}
As in Lemma 3.5 of \cite{jo-kribs} we have 
\begin{equation*}
\left( S_{i-1}^{*}S_{j-1}\right) =\left( AA^{*}\right) _{i,j}
\end{equation*}
and similarly 
\begin{equation*}
\left( \tilde S_{i-1}^{*}\tilde S_{j-1}\right) =\left( \tilde{A}\tilde{A}%
^{*}\right) _{i,j}
\end{equation*}

For $1\le i\le N$, let $T_i$ be the $i$-th creation operators on the twisted
Fock space.

To find the action of $\tilde T_i^{*}$ on the spanning vectors we consider 
\begin{eqnarray*}
\langle \tilde T_i ( w\otimes h ) \mid \tilde T_j ( w^{\prime }\otimes
h^{\prime } ) \rangle &=&\langle iw\otimes h\mid \left( jw^{\prime }\right)
\otimes h^{\prime }\rangle \\
&=& \left\langle h \biggm| \sum_\sigma q^{2\sigma \left( i\right) }\Phi
\left( e_{i_1,i_1^{\prime }}\cdots e_{i_k,i_k^{\prime }}\right) h^{\prime }
\right\rangle
\end{eqnarray*}
implying 
\begin{eqnarray*}
\tilde T_j^{*}\tilde T_i &=&\delta _{i,j}\sum_{\sigma \in S_N}q^{2\sigma
\left( i\right) } \\
\ &=&\frac{(1-q^{2N})}{1-q^2}\delta _{i,j}
\end{eqnarray*}
We consider $S_i$ and $\tilde S_i$ coming from wavelet analysis.

For $1\le i\le N$ let $\tilde T_i$ and $\tilde S_i$ be the creation
operators and the operator coming from the wavelet construction.

Then 
\begin{equation*}
\tilde S_{i-1}^{*}\tilde S_{j-1} =\left( \tilde{A}\tilde{A}^{*}\right)
_{i,j}= \tilde{T}_i\tilde{T}_j^{*}
\end{equation*}

Let us describe the action of $\tilde{T}_{i}T_{j}^{\ast }$ on the spanning
vectors. Since 
\begin{equation*}
\langle T_{i}^{\ast }(w\otimes h)\mid w^{\prime }\otimes h^{\prime }\rangle
=\langle w\otimes h\mid \left( iw^{\prime }\right) \otimes h^{\prime }\rangle
\end{equation*}%
Thus $\langle \left( \smash{\tilde{T}_iT_j^{*}}\right) \left( w\otimes
h\right) \mid w^{\prime }\otimes h^{\prime }\rangle $ can be computed as
follows: 
\begin{eqnarray*}
&&\langle \left( \smash{\tilde{T}_iT_j^{*}}\right) \left( w\otimes h\right)
\mid w^{\prime }\otimes h^{\prime }\rangle \\
&&\qquad =\langle \langle i,j\rangle ^{-1}w\otimes h\mid w^{\prime }\otimes
h^{\prime }\rangle \\
&&\qquad \qquad \langle p_{i,j}^{-1}w_{i_{1},i_{1}^{\prime
}},w_{i_{2},i_{2}^{\prime }},\dots ,w_{i_{k},i_{k}^{\prime }}\otimes h\mid
w_{i_{1},i_{1}^{\prime }}^{\prime },w_{i_{2},i_{2}^{\prime }}^{\prime
},\dots ,w_{i_{k},i_{k}^{\prime }}^{\prime }\otimes h^{\prime }\rangle \\
&&\qquad =\left\langle \frac{\left( 1-q^{-N}\ \right) }{\left(
1-q^{-1}\right) }\left( 1-q^{-1}\right) w_{i_{1},i_{1}^{\prime
}},w_{i_{2},i_{2}^{\prime }},\dots ,w_{i_{k},i_{k}^{\prime }}\otimes h\biggm|%
\right. \\
&&\qquad \qquad \qquad \qquad \qquad \qquad \qquad \qquad \left. \biggm|%
w_{i_{1},i_{1}^{\prime }}^{\prime },w_{i_{2},i_{2}^{\prime }}^{\prime
},\dots ,w_{i_{k},i_{k}^{\prime }}^{\prime }\otimes h^{\prime }\right\rangle
\end{eqnarray*}%
On the other side we have 
\begin{eqnarray*}
&&\langle T_{j}^{\ast }\tilde{T}_{i}\left( w\otimes h\right) \mid w^{\prime
}\otimes h^{\prime }\rangle \\
&&\qquad =\langle \left( iw\right) \otimes h\mid \left( jw^{\prime }\right)
\otimes h^{\prime }\rangle \\
&&\qquad =\left\langle \frac{\left( 1-q^{N}\ \right) }{\left( 1-q\right) }%
\left( 1-q\right) w_{i_{1}},w_{i_{2}},\dots ,w_{i_{k}}\otimes h\biggm|%
w_{i_{1}^{\prime }}^{\prime },w_{i_{2}^{\prime }}^{\prime },\dots
,w_{i_{k}^{\prime }}^{\prime }\otimes h^{\prime }\right\rangle
\end{eqnarray*}%
thus it follows: 
\begin{equation*}
T_{j}^{\ast }\tilde{T}_{i}=1-q^{N}
\end{equation*}%
and 
\begin{equation*}
\tilde{T}_{i}T_{j}^{\ast }=1-q^{-N}
\end{equation*}%
Hence: 
\begin{equation*}
\tilde{T}_{i}T_{j}^{\ast }-T_{j}^{\ast }\tilde{T}_{i}=\delta _{i,j}[N]_{q}%
\mathbf{1}
\end{equation*}

Then the creation and annihilation operators $\tilde{T}_{i}$ and $%
T_{j}^{\ast }$ and $N$, viewed as a number operator, yield a representation
of the $q$-oscillator algebra.
\end{proof}

\noindent \textbf{Concluding remarks:}

We have shown how tools from transform theory and wavelet analysis help us
in the construction of new representations of certain $q$-relations
fromstatistical mechanics. The $q$-relations have been studied earlier in
connection with quantum fields \cite{Bo-Speicher} and statistical mechanics 
\cite{RW}. In particular, the paper \cite{RW} serves to show that the $q$%
-relations interpolate between the bosons and the fermions. Further, in \cite%
{RW}, the partition function is calculated for the quons and it is
established that it exhibits Gibbs' paradox. As a result, the corresponding
notions of entropy, free energy and particle number break with our
traditional understanding of thermodynamical quantities.


\frenchspacing

\begin{thebibliography}{99}
\bibitem{Bo-Speicher} M. Bo\.zejko, R. Speicher, An example of a generalized 
{B}rownian motion, \emph{Comm. Math. Phys.} \textbf{137} (1991) 519--531.

\bibitem{jo-bra-ev} O. Bratteli, D.E. Evans, P.E.T. Jorgensen, Compactly
supported wave\-lets and representations of the Cuntz relations, \emph{Appl.
Comput. Harmon. Anal.} \textbf{8} (2000), 166--196.

\bibitem{jo-bra} O. Bratteli, P.E.T. Jorgensen, Isometries, shifts, Cuntz
algebras and multiresolution wavelet analysis of scale $N$, \emph{Integral
Equations Operator Theory} \textbf{28} (1997), 382--443.

\bibitem{pr-vy} Vyjayanthi Chari, Andrew N. Pressley, \emph{A Guide to
Quantum Groups}, Cambridge University Press, Cambridge, 1995, corrected
reprint of the 1994 original.

\bibitem{dau} I. Daubechies, \emph{Ten Lectures on Wavelets}, CBMS-NSF
Regional Conf. Ser. in Appl. Math., vol. 61, Society for Industrial and
Applied Mathematics, Philadelphia, 1992.

\bibitem{jo-kribs} P.E.T. Jorgensen, D. Kribs, Wavelet representations and
Fock space on positive matrices, http://arXiv.org/abs/math.CA/0204034, to
appear in \emph{J.~Funct. Anal.}

\bibitem{jo-pa} P.E.T. Jorgensen, A. Paolucci, Multiresolution wavelet
analysis of Bessel functions of scale $\nu +1$,
http://arXiv.org/abs/math.FA/0006103

\bibitem{jo-we} P.E.T. Jorgensen, L.M. Schmitt, R.F. Werner, $q$-canonical
commutation relations and stability of the Cuntz algebra, \emph{Pacific J.
Math.} \textbf{165} (1994), 131--151.

\bibitem{BrJo02b} O.~Bratteli and P.E.T. Jorgensen, \emph{Wavelets through a 
{L}ooking {G}lass: The {W}orld of the {S}pectrum}, Applied and Numerical
Harmonic Analysis, Birkh\"auser, Boston, 2002.

\bibitem{IS} I. Sneddon, \emph{Fourier Transforms}, Dover, New York, 1995,
reprint of the 1951 original, McGraw-Hill, New York.

\bibitem{RS} R.F. Swarttouw, \emph{The Hahn-Exton $q$-Bessel function},
Ph.D. thesis, Delft University of Technology, 1992 (see also R. F.
Swarttouw, H. G. Meijer, A $q$-analogue of the Wronskian and a second
solution of the Hahn-Exton $q$-Bessel difference equation, \emph{Proc. Amer.
Math. Soc.} \textbf{120} (1994), 855--864).

\bibitem{RW} R.F. Werner, The free quon gas suffers Gibbs paradox, Phys.
Rev. D. \textbf{48} (1993), 2929--2934.
\end{thebibliography}
\end{document}